\documentclass[conference]{IEEEtran}
\IEEEoverridecommandlockouts
\usepackage{graphicx}
\usepackage{amsmath}
\usepackage{amssymb}
\usepackage{algorithmicx}
\usepackage[ruled]{algorithm}
\usepackage{algpseudocode}
\usepackage{setspace}

\def\C{\mathbb{C}}
\newcommand{\supp}{\mathrm{supp}}
\newcommand{\sgn}{\mathrm{sign}}
\IEEEaftertitletext{\vspace*{-5mm}}
\begin{document}
%
\title{One-Bit Compressive Sensing with Partial Support}

\author{\IEEEauthorblockN{Phillip North}
\IEEEauthorblockA{Dept. of Mathematical Sciences\\
Claremont McKenna College\\
Claremont CA 91711\\
Email: PNorth15@students.claremontmckenna.edu}
\and
\IEEEauthorblockN{Deanna Needell}
\IEEEauthorblockA{Dept. of Mathematical Sciences\\
Claremont McKenna College\\
Claremont CA 91711\\
Email: dneedell@cmc.edu}\thanks{D.N.\ was partially supported by NSF CAREER DMS-1348721 and the Alfred P.\ Sloan Foundation. }
}

\maketitle

\begin{abstract}
The Compressive Sensing framework maintains relevance even when the available measurements are subject to extreme quantization, as is exemplified by the so-called one-bit compressed sensing framework which aims to recover a signal from measurements reduced to only their sign-bit.  In applications, it is often the case that we have some knowledge of the structure of the signal beforehand, and thus would like to leverage it to attain more accurate and efficient recovery.  This work explores avenues for incorporating such partial-support information into the one-bit setting.  
Experimental results demonstrate that newly proposed methods of this work yield improved signal recovery even for varying levels of accuracy in the prior information. This work is thus the first to provide recovery mechanisms that efficiently use prior signal information in the one-bit reconstruction setting.
\end{abstract}

\IEEEpeerreviewmaketitle

\section{Introduction}
\textit{Compressed Sensing} (CS) addresses the problem of accurately acquiring high dimensional signals from a set of relatively few linear measurements \cite{RefWorks:59} \cite{RefWorks:113} \cite{RefWorks:68}.  The problem can be formulated mathematically via the system $\tilde{y} = \Phi x$, where $\Phi \in \C^{m \times n}$ is the \textit{measurement matrix}.  In the compressed setting, $m\ll n$ but one utilizes the assumption that the signal $x$ possesses some additional structure, such as \textit{sparsity}; we say that $x$ is \em $k$-sparse\em \ when
$$
\|x\|_0 := |\supp(x)| = k \ll n.
$$
CS has seen a vast amount of progress (see e.g. \cite{eldar2012compressed,foucart2013}), and it is now well-known that for suitable matrices $\Phi$ (for example i.i.d. Gaussian), any $k$-sparse vector ${x}$ can be recovered from $\tilde{y}\in\C^m$ when $m \approx k\log(n/k)$.

Unfortunately, the majority of theoretical work in CS assumes that the measurements are acquired with infinite precision whereas in practice they must be quantized.  The extreme quantization setting where only the sign bit is acquired is known as \textit{one-bit compressed sensing} \cite{Boufounos2008}.  In this framework, the measurements now take the form ${y}_i = \sgn(\langle x, \phi_i \rangle)$ where $\phi_i$ denotes the $i$th row of the measurement matrix $\Phi$.  Typically one then loses the ability to recover the magnitude of $x$ and thus assumes the signal has a fixed norm (e.g. unit-norm), although there are adaptive techniques to overcome this as well \cite{RefWorks:348,exponentialBFNPW14}.

\subsection{Existing one-bit methods}
Although one-bit CS is a relatively new technology, efficient recovery algorithms have been studied.  There are mainly two types of methods, those based on linear programming \cite{pv-1-bit,pv-noisy-1bit,gnjn2013} and on iterative approaches \cite{jacques2013quantized,biht}.  In this work we focus on the iterative approach, and in particular the Binary Iterative Hard Thresholding Algorithm (BIHT) \cite{biht}, an extension of the traditional Iterative Hard Thresholding Algorithm (IHT) \cite{PaperIHT}. 
Assuming that our desired signal is $k$-sparse, the objective of BIHT is to return a solution that is $k$-sparse and consistent with the given sign measurements.  Viewed as solving an optimization problem, at each iteration BIHT computes and takes a step in the direction of the gradient to attain a new approximation. This approximation is then thresholded to retain only the $k$ largest in magnitude entries. Finally, after a consistent approximation is attained or enough iterations have elapsed, the estimation is normalized and returned. Algorithm \ref{alg:BIHT} presents a more detailed explanation of BIHT.  Here and throughout we use the thresholding function $prune(z, k)$ which returns the vector $z$ with all but the $k$ largest in magnitude entries set to zero.

\begin{algorithm}[ht]
    \caption{Binary Iterative Hard Thresholding (BIHT). Given: measurement matrix $\Phi$, one-bit measurements ${y}$, assumed sparsity level $k$, gradient step-size $\tau$}\label{alg:BIHT}
\begin{algorithmic}[1]
\Procedure{BIHT}{$\Phi$, ${y}$, $k$} 
\State $\tilde{x} = 0$ \Comment{Initialize trivial approximation of $x$}
\Repeat
    \State $\Gamma = \tilde{x} + \frac{\tau}{2} \Phi' ({y} - sign(\Phi x_i)) $ \Comment{Gradient step} 
    \State $\tilde{x} = prune(\Gamma, k)$ \Comment{Hard threshold} 
\Until halting criterion satisfied
\State \Return     $\frac{\tilde{x}}{\| \tilde{x} \|_2}$ \Comment{Normalize}
\EndProcedure
\end{algorithmic}
\end{algorithm}

\section{Methods using partial support}

In many applications, it is not only known that the signal of interest is sparse, but additional information about the support of the signal may also be known.  For example, it is well-known that the support of wavelet representations of natural images largely resides in the low frequencies.  In distributed settings, the signals of interest may be highly correlated so that partial support information may be obtained from neighboring atoms.  This framework can be modeled by assuming that there is some partial estimate $\tilde{T}$ of the signal support $T := \supp(x)$ a priori.  The work of Mansour et.al. \cite{DBLP:journals/corr/MansourS14} \cite{DBLP:journals/corr/abs-1010-4612} extend the conventional $\ell_1$-minimization method to a \textit{weighted} $\ell_1$ approach that effectively incorporates such a support estimate.  
To incorporate the estimate $\tilde{T}$, consider the following program:
$$
\min_x \sum_{i=1}^n w_i | x_i | \;\;\text{s.t.}\;\; \tilde{y} = \Phi x 
\text{   with     } w_i = \begin{cases} c \in \left[0,1\right] \text{ if } w_i \in \tilde{T} \\ 1 \text{ if } w_i \notin \tilde{T}. \end{cases}
$$
This program now imposes a penalty for placing non-zero entries in locations not specified in the support estimate $\tilde{T}$. The value $c$ in the weight vector $w$ can be determined by the confidence in each element of $\tilde{T}$.  We build upon this notion of weighting the estimation vector by developing analogous methods for the one-bit setting, focusing on iterative methods.

\subsection{Oracle estimation}
\begin{figure}[ht]
\centering
\includegraphics[scale = .25]{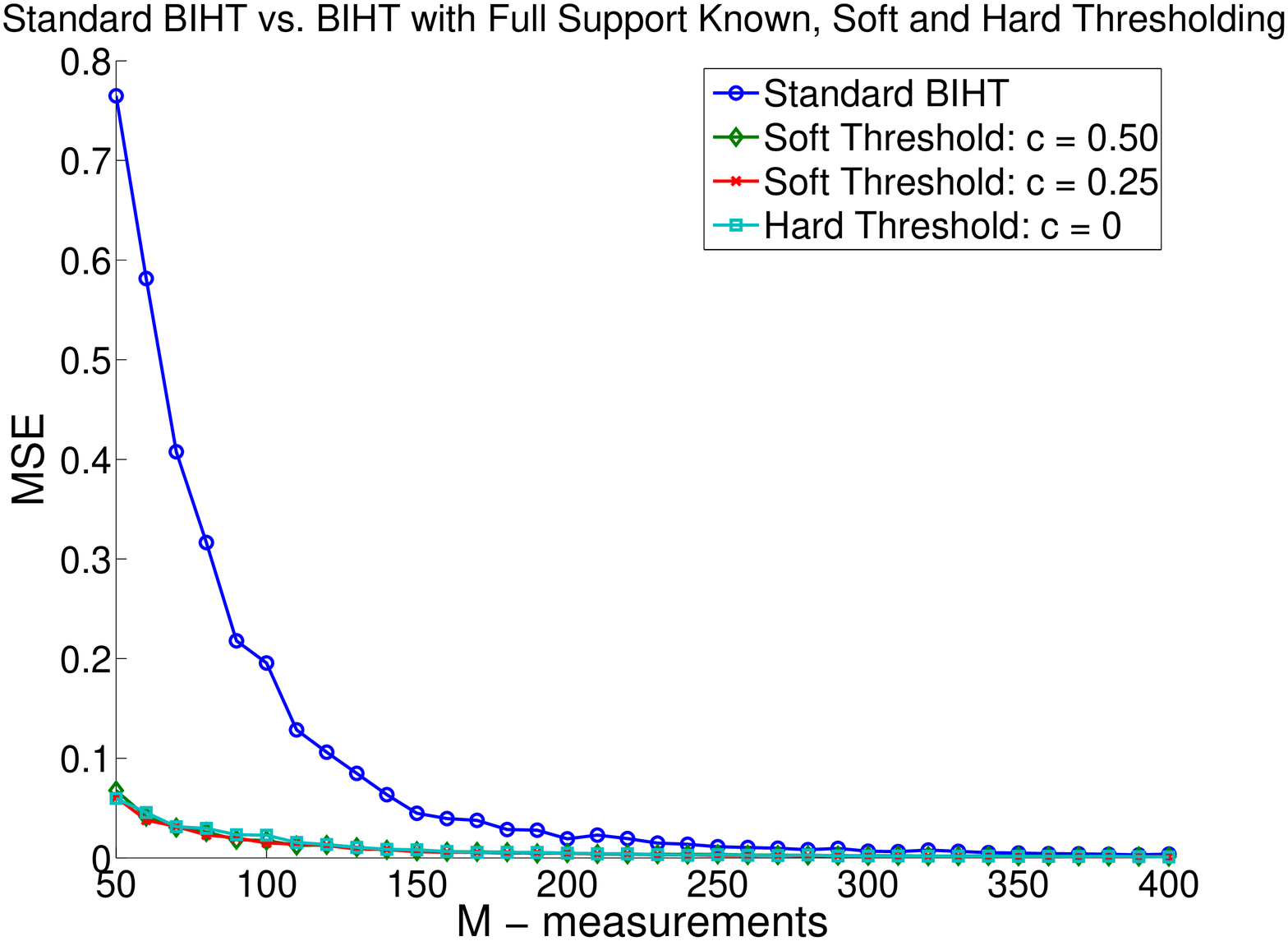}
\caption{Performance of BIHT when support estimate is exact and we employ soft and hard thresholding ($c=0$).}\label{fig:oracle}
\end{figure}

As a first, most basic approach, let us assume that our support estimate $\tilde{T}$ is completely accurate, i.e. $\tilde{T} = T$. At the pruning step, no matter our result, we could simply set all entries of $\Gamma$ not in $\tilde{T}$ to zero; instead of locating the $k$ largest entries of $\Gamma$, we would naively only retain the entries of $\tilde{T}$. The entries of our new estimate would be determined as follows:
$$
\tilde{x}_i = \begin{cases} \Gamma_i \text{ if } i \in \tilde{T} \\ 0 \text{ if } i \in \tilde{T}^C. \end{cases} 
$$
 As a similar approach, we could instead try soft-thresholding the entries of $\Gamma$ not in $\tilde{T}$, multiplying the entries of $\Gamma$ not in $\tilde{T}$ by some constant $0<c<1$. Figure \ref{fig:oracle} shows the result of this approach; here and throughout, unless otherwise noted the signal length is $n=256$, sparsity level $k = 8$, $\Phi$ has standard normal entries, $\tau = 0.001$, the support of the signal $x$ is distributed uniformly at random, the magnitudes of the non-zero entries are standard normal, and the algorithm is run until the estimate changes by less than $10^{-10}$ or after 1000 iterations.  Figure \ref{fig:oracle} shows the mean squared error (MSE) averaged over 100 trials for various values of $m$, the number of measurements.  At each iteration, entries of $\Gamma$ that are not in the support estimate are scaled down by a factor $c$, so that when $c = 0$ the result is simple hard-thresholding.  Observe unsurprisingly how powerful a perfect support estimate can be via this bold hard-thresholding strategy.  Figure \ref{fig:oracle} also shows that soft-thresholding with the perfect support estimate performs identically to hard-thresholding. This result is not hard to understand: after $j$ iterations the elements not in $\tilde{T}$ have been scaled down by $c^j$, which approaches zero after many iterations. Of course, knowing the full support beforehand is not likely, but these examples show that this information can seriously expedite the recovery of $x$.

\subsection{Soft thresholding}
\begin{figure}[ht]
                (a) \includegraphics[scale = .25]{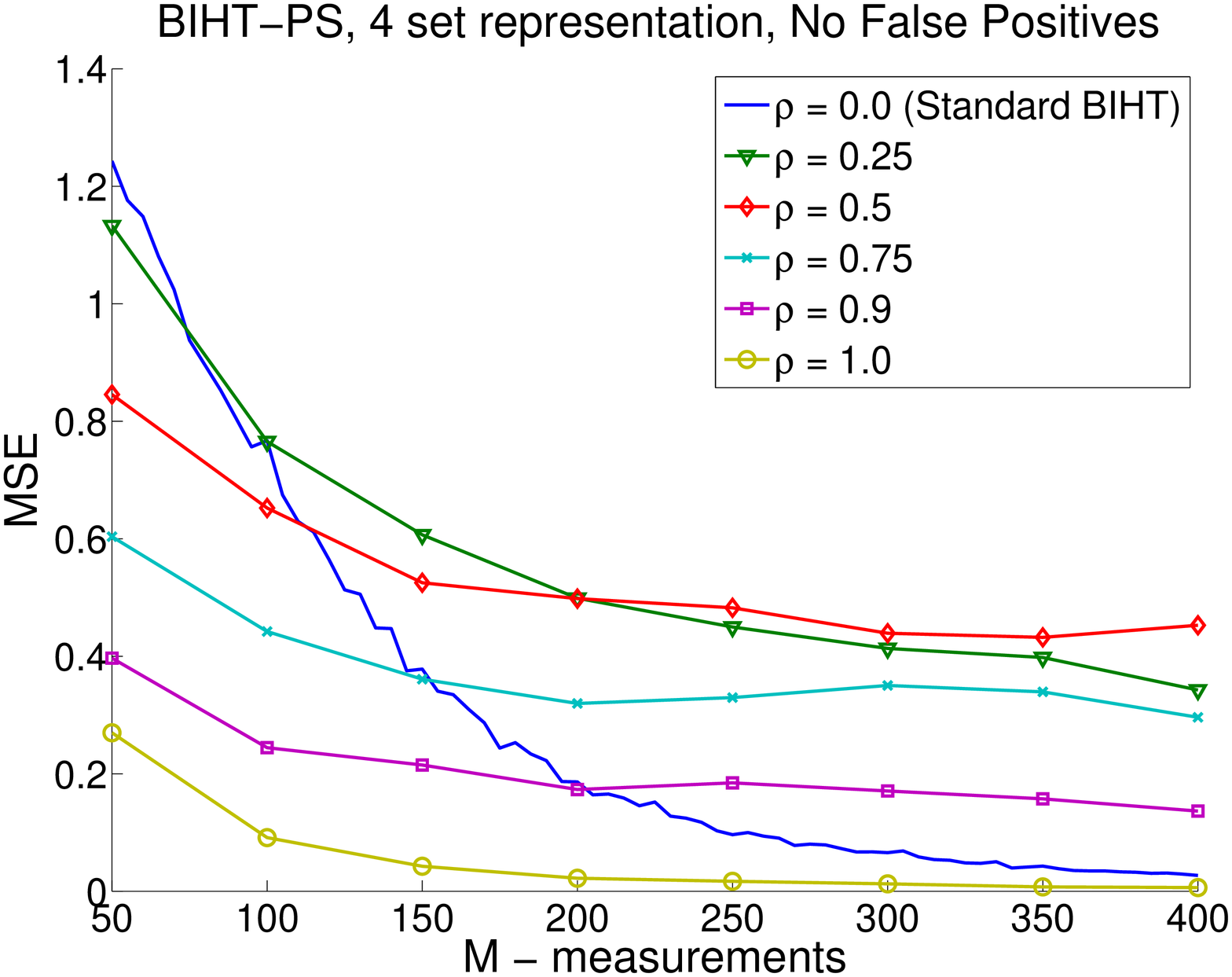}\\
                (b)\includegraphics[scale = .25]{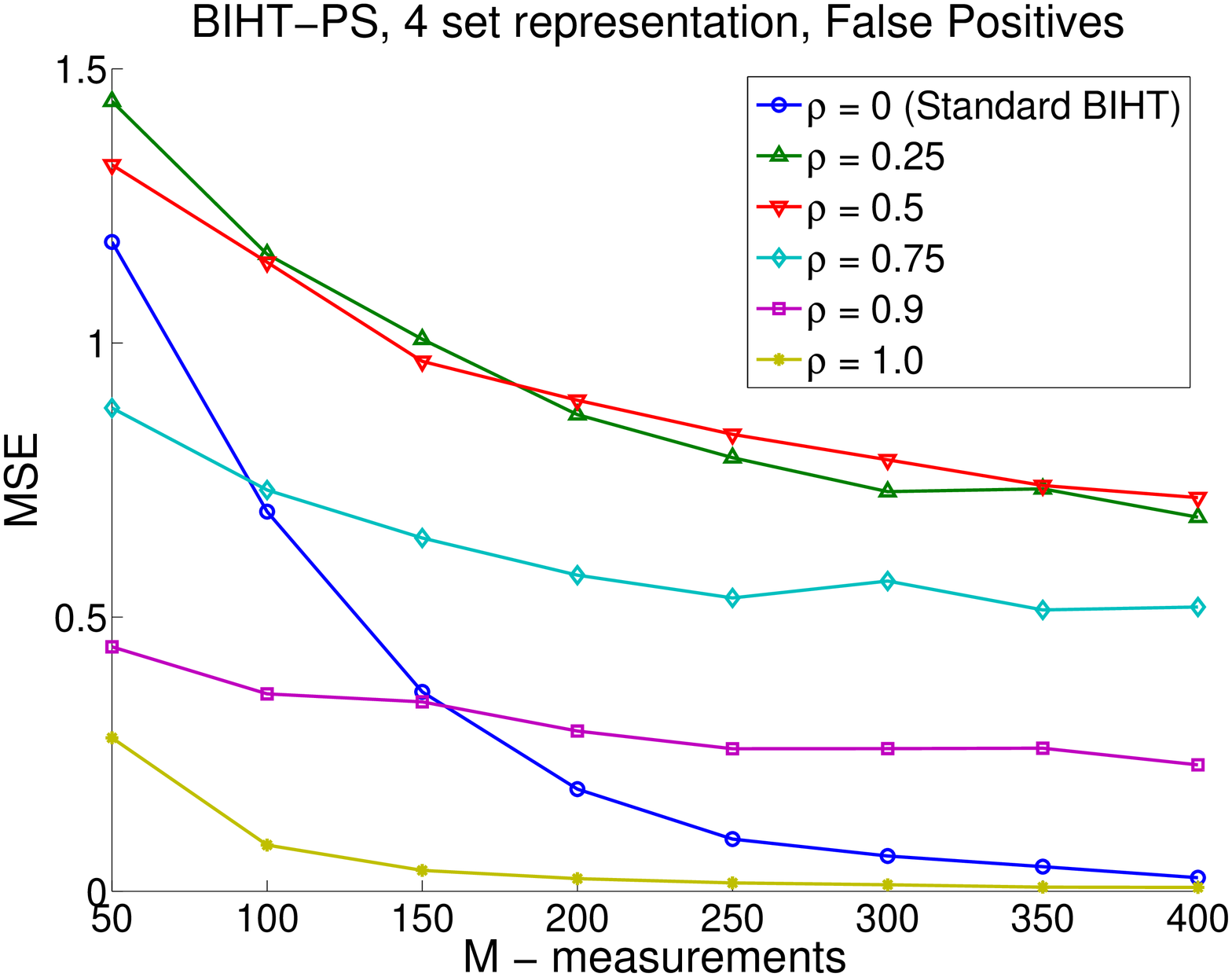}
        \caption{Performance of BIHT using soft thresholding. In (a) we know the estimate to contain $\rho k$ correct elements (no false positives), in (b) there are an additional $(1-\rho)k$ incorrect elements. }\label{fig:soft}
\end{figure}

In practice our support estimate will seldom equal the true support, and thus at a given iteration we should consider both $\tilde{T}$ and the locations of the $k$ greatest in magnitude entries of $\Gamma$, denoted $\dot{T}$, when updating our estimate. Let us assume that we know the support estimate to be $(\rho \times 100 )\%$ accurate. Now, in the pruning step of BIHT the entries of $\Gamma$ may be thresholded according to four distinct possibilities:
$$
\tilde{x}_i = \Gamma_i w_i \text{    where    } w_i = \begin{cases} 
1 \text{ if } i \in \tilde{T} \cap \dot{T} \\
1 \text{ if } i \in \tilde{T} \cap \dot{T}^C \\
1-\rho \text{ if } i \in \tilde{T}^C \cap \dot{T} \\
0 \text{ if } i \in \tilde{T}^C \cap \dot{T}^C. \\
\end{cases} 
$$
Additionally, we can use this 4-set framework when our support estimate includes erroneous elements. Suppose we have a support estimate that contains $(\rho \times k)$ correct elements but also includes $(1-\rho)\times k$ incorrect elements. Figure \ref{fig:soft} shows BIHT's performance when using this 4-set representation to incorporate prior support information. In the case of no false positives, when no elements of our support estimate are incorrect, we see that performance of BIHT is not bad: as $\rho$ increases the MSE decays. Similarly, with the inclusion of false positives, performance is intuitive: improvements are seen when there are more correct estimates than incorrect estimates, i.e. $\rho \geq 0.5$; for lesser values of $\rho$ the support estimate consists mostly of incorrect elements and performance is worse than standard BIHT. Also note that under this 4-set representation a $k$-sparse approximation is not necessarily returned. If $\tilde{T} \cap \dot{T} = \emptyset$, then in fact a $2k$-sparse solution is returned, certainly this will result in a less accurate solution. Perhaps this is the reason for the slower rate of error decay in comparison with standard BIHT, as seen in Figure \ref{fig:soft}. 

\subsection{Supervised weighting}
\begin{figure}[ht]
                (a)\includegraphics[scale = .25]{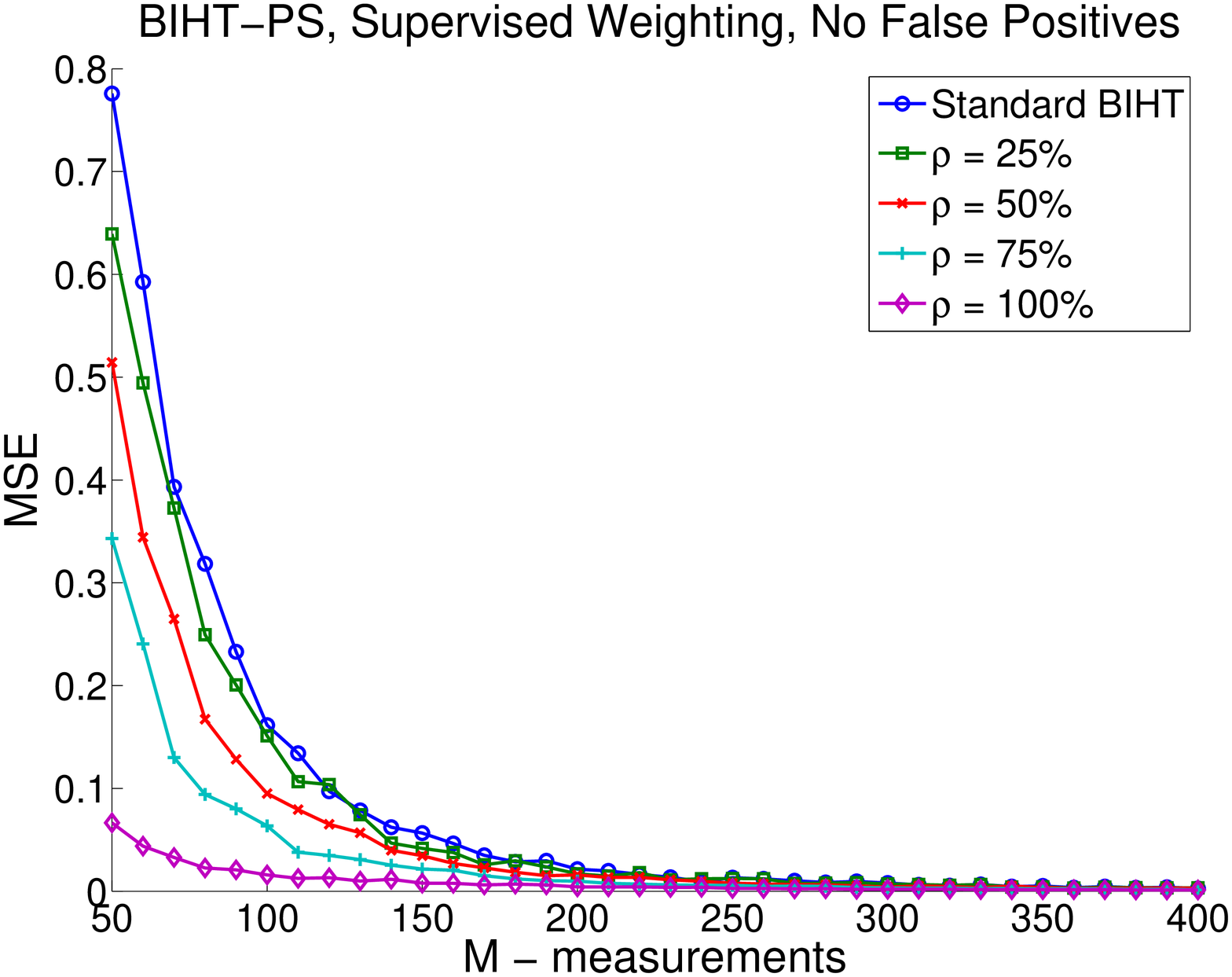}\\
                (b)\includegraphics[scale = .25]{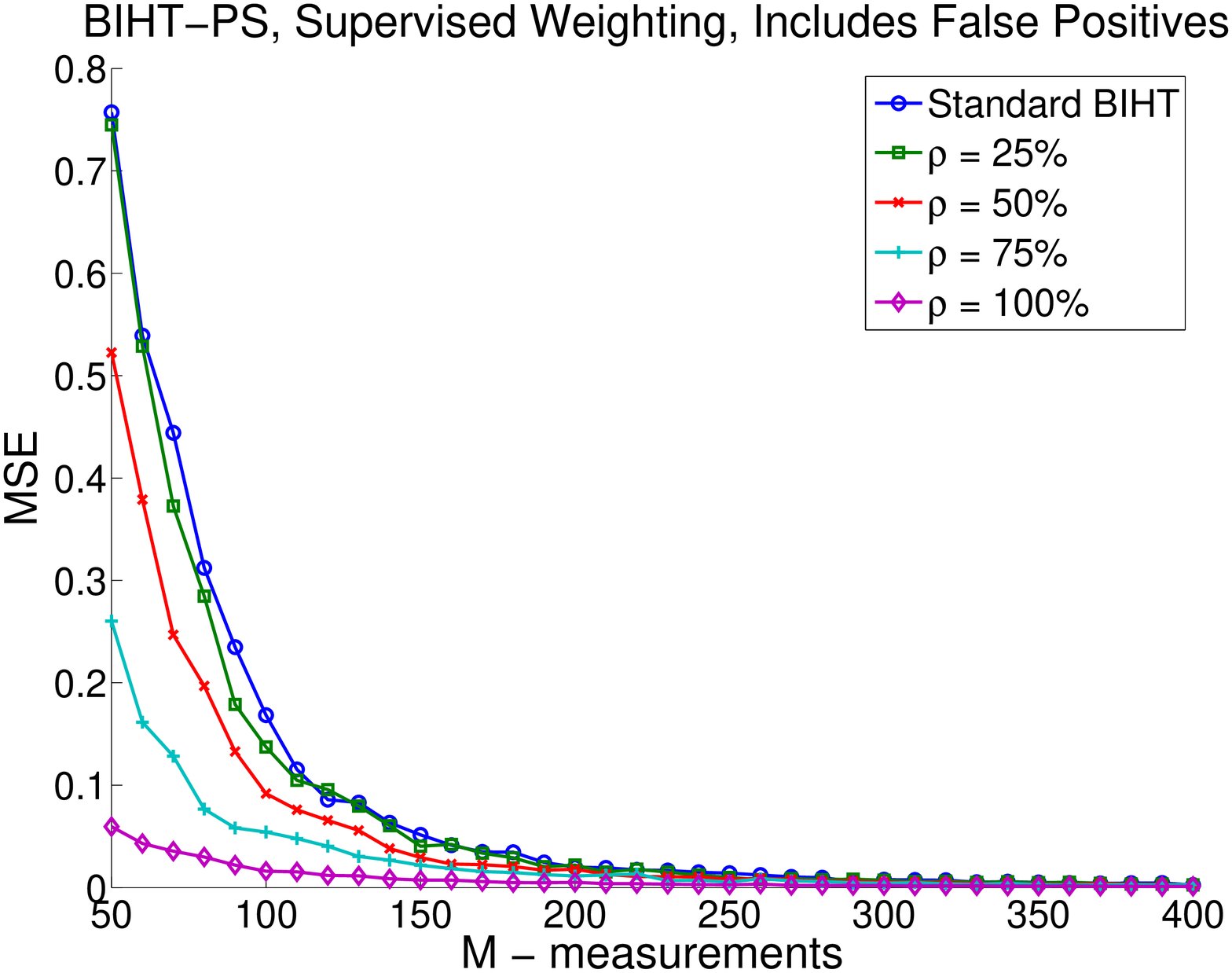}
        \caption{Performance of BIHT when using supervised weighting. In (a) we know the estimate to contain $\rho k$ correct elements (no false positives), in (b) there are an additional $(1-\rho)k$ incorrect elements. }\label{fig:super1}
\end{figure}

As a means for incorporating a partial support estimate and returning a $k$-sparse solution, we may use the weighting framework for traditional $\ell_1$-minimization as presented in \cite{DBLP:journals/corr/MansourS14}.  We refer to this model as \textit{supervised}, since the partial support estimate $\tilde{T}$ must be obtained a priori, external to the algorithm itself.  Again, suppose we believe our estimate $\tilde{T}$ to be $(\rho \times 100)\%$ accurate. Let us create a weight vector $w$ where 
$$
w_i = \begin{cases} 1 \text{ if } i \in \tilde{T} \\ 1-\rho \text{ if } i \in \tilde{T}^C.  \end{cases}
$$
Then, at some iteration of BIHT, once we compute $\Gamma$ let us multiply component-wise $\Gamma$ and $w$ to attain $\psi = \Gamma \odot w$. Now we locate the $k$ largest elements of $\psi$, denoted $T_{\psi} := prune(\psi, k)$, and hard threshold the elements of $\Gamma \in T_{\psi}^C$, i.e. set them equal to zero. This approach is very similar to BIHT when there is no prior support estimate, except that here $\Gamma$ is instead pruned to retain the $k$ greatest in magnitude elements of $\Gamma \odot w$. A more thorough break down of this procedure is presented in Algorithm \ref{alg:super}. We see that this supervised weighting approach outperforms the 4-set soft thresholding formulation both when false positives are and are not included. Figure \ref{fig:super1} shows that for every value of $m$, when some prior information exists, \textit{the weighting approach performs the same as or better than standard BIHT.} 

\begin{figure}[ht]
\centering
\includegraphics[scale = .25]{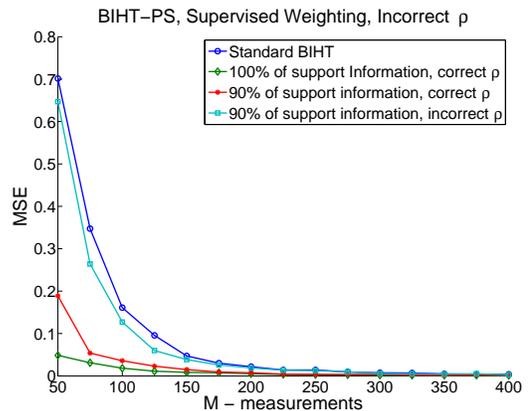}
\caption{Effect of inaccurate choice of $\rho$ in supervised weighting.  For incorrect $\rho$ (cyan), $\rho=0.1$ (rather than $\rho=0.9$) was used in the weight vector.}\label{fig:super2}
\end{figure}

\begin{algorithm}[ht]
    \caption{Binary Iterative Hard Thresholding with Partial Support Estimate Weighting (BIHT-PSW). Given: measurement matrix $\Phi$, one-bit measurements ${y}$, sparsity level $k$, support estimate $\tilde{T}$, accuracy of support estimate $\rho$, step-size $\tau$}\label{alg:super}
\begin{algorithmic}[1]
\Procedure{BIHT-PSW}{$\Phi$, ${y}$, $k$, $\tilde{T}$, $\rho$} 
\State $w_i = \begin{cases} (1-\rho) \text{   if   } i \in \tilde{T}^C \\ 1 \text{   if   } i \in \tilde{T} \end{cases}$ \Comment{Construct weights}
\State $\tilde{x} = 0$ \Comment{Initialize trivial approximation of $x$}
\Repeat
    \State $\Gamma = \tilde{x} + \frac{\tau}{2} \Phi' ({y} - sign(\Phi x^i)) $ \Comment{Gradient step} 
    \State $\Omega = prune(\Gamma \odot w, k)$ \Comment{Prune weighted update}
    \State $\tilde{x}_j = \begin{cases} \Gamma_j \text{   if   } j \in \Omega \\ 0 \text{   if   } j \in \Omega^C \end{cases}$ \Comment{Update approximation}
\Until halting criterion is satisfied
\State \Return     $\frac{\tilde{x}}{\| \tilde{x} \|_2}$ \Comment{Normalize}
\EndProcedure
\end{algorithmic}
\end{algorithm}
 
This weighting framework for incorporating partial support information into BIHT (BIHT-PSW) performs well in the current context. However, we are assuming that the value for $\rho$ is correct. This is a very bold assumption that, in practice, is not likely to hold. In the weighting step, the value of $\rho$ determines how diminished the magnitude of an entry off the support estimate will be. If we are more confident in certain entries being non-zero, then other entries will be scaled by a constant closer to zero. Empirical results, as displayed in Figure \ref{fig:super2}, show that using the correct value of $\rho$ is crucial to the incorporation of a partial support estimate.

\subsection{Unsupervised re-weighting}
\begin{figure}[ht]
               (a) \includegraphics[scale = .25]{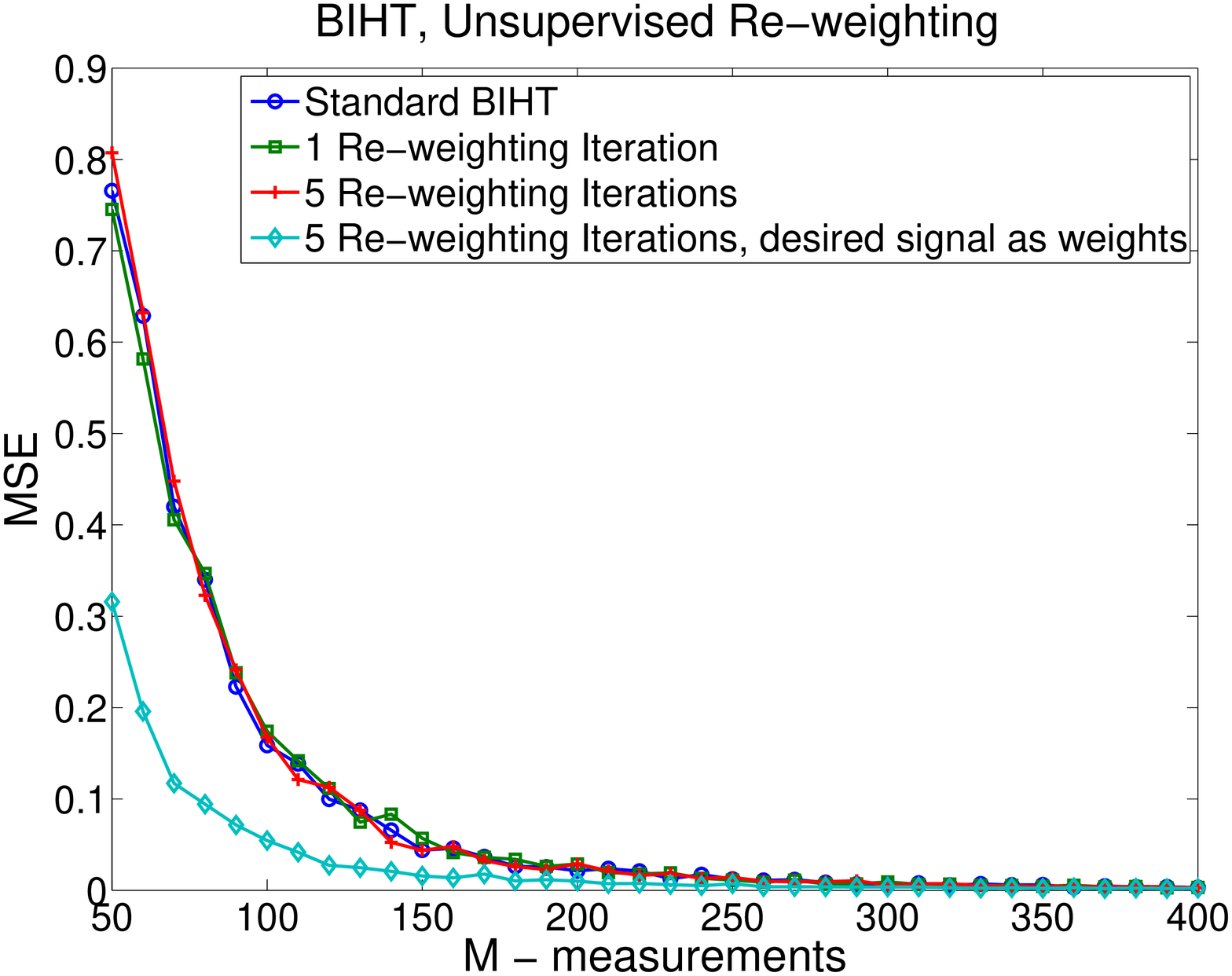}\\
               (b) \includegraphics[scale = .25]{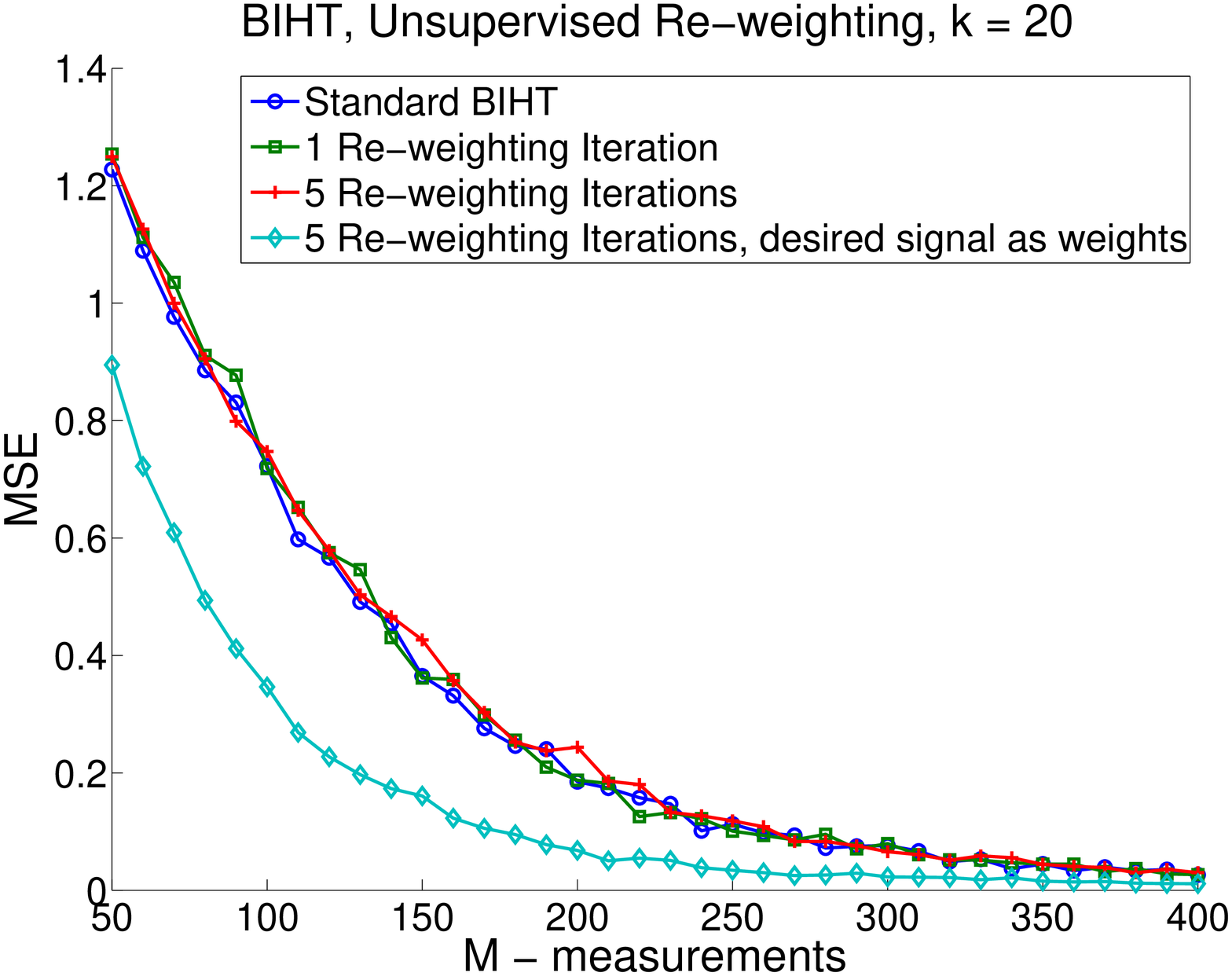}
        \caption{Performance of BIHT-URW for various re-weighting iterations. In (a) the sparsity level is $k=8$, and in (b) $k=20$. In both figures, the cyan ($\lozenge$) line is a result of creating a weight vector out of the desired signal $x$. }\label{fig:unsuper}
\end{figure}

  The improvements of this approach prompt one to ask whether such a result can be leveraged even when no support estimate is available. If we take $m$ measurements of $x$ and use BIHT to attain some estimate $\tilde{x}$, then we may use the support of $\tilde{x}$, denoted $\ddot{T}$, as an estimate for $T$. Then, we could use $\ddot{T}$ to run BIHT-PSW, getting a more accurate approximation, and so on. This is reminiscent of the re-weighted $\ell_1$-minimization approach in classical compressed sensing \cite{RefWorks:308}, and is the general idea behind the BIHT Unsupervised Re-weighting (BIHT-URW) algorithm, presented as Algorithm \ref{alg:unsuper}.  We refer to this model as \textit{unsupervised}, since no outside information about the support is required.  In this case, since $\rho$ is certainly unknown, we utilize a parameter $\lambda$ in place of $\rho$.  The performance of BIHT-URW is displayed in Figure \ref{fig:unsuper}.  Unfortunately, we observe no improvement from standard BIHT; we conjecture this is perhaps due to the arbitrary selection of $\lambda$ (and thus $\rho$) within the method, which we have set to $\lambda=0.5$ in Figure \ref{fig:unsuper}.  As a benchmark for optimal performance of BIHT-URW we tested the method using a weight vector from the actual signal $x$ (rather than its approximation), and then ran the algorithm as usual. This resulted in a significant improvement in performance, as shown in the cyan curve of Figure \ref{fig:unsuper}.  Of course, this method is not applicable in practice, but demonstrates the potential of such an approach if better ways of estimating $\lambda$ can be obtained, possibly adapting from iteration to iteration.

\begin{algorithm}[ht]
    \caption{Binary Iterative Hard Thresholding with Unsupervised Re-weighting. Given: measurement matrix $\Phi$, one-bit measurements ${y}$, sparsity level $k$, step-size $\tau$, accuracy parameter $\lambda$, number of re-weighting iterations $n$}\label{alg:BIHT-URW}\label{alg:unsuper}
\begin{algorithmic}[1]
\Procedure{BIHT-URW}{$\Phi$, ${y}$, $k$, $n$} 
\State $\ddot{T} = \supp(BIHT(\Phi, {y}, k))$ \Comment{BIHT support estimate} 
\Repeat 
    \State $\tilde{x} = \text{BIHT-PSW}(\Phi, {y}, k, \ddot{T}, \lambda)$ \Comment{Estimate} 
    \State $\ddot{T} = \supp(\tilde{x})$ \Comment{Update support estimate} 
\Until $n$ iterations have completed
\State \Return     $\frac{\tilde{x}}{\| \tilde{x} \|_2}$ \Comment{Normalize} 
\EndProcedure
\end{algorithmic}
\end{algorithm}

\section{Conclusion}

We presented several weighting variants of the BIHT algorithm for one-bit CS when partial support estimation is known.  We demonstrate that our methods effectively utilize the support information, but that leveraging support estimates in an unsupervised fashion is not straightforward.  We believe future work in this area could lead to re-weighted methods which improve upon existing approaches.

\bibliographystyle{ieeetran}
\bibliography{bib}

\end{document}